\documentclass[20pt]{article}
 \usepackage[T2A]{fontenc}
\usepackage[margin=1in]{geometry}
\usepackage{dsfont}
\usepackage{amsmath, amsthm, amssymb, mathrsfs}
\usepackage[abbrev,non-sorted-cites]{amsrefs}
\usepackage[all]{xy}
\usepackage{graphicx}
\usepackage{extarrows}
\usepackage{hyperref}
\usepackage{xcolor}
\usepackage{tikz-cd}
\renewcommand{\thefootnote}{\arabic{footnote}}

\begin{document}

\begin{Large}
\centerline{A note on general isolation result in Diophantine Approximation}
\end{Large}
\vskip+1cm

\begin{Large}
\centerline{ 
\begingroup
\renewcommand{\thefootnote}{\raisebox{0.5ex}{\small$\dagger$}}
Sergei Pitcyn\footnote{
Moscow State University and  Moscow Center of Fundamental and Applied Mathematics,
picyn98@mail.ru 
}
\endgroup 
\begingroup
\renewcommand{\thefootnote}{\raisebox{0.5ex}{\small$\ddagger$}}
 and Nikolay Moshchevitin\footnote{
 Institut für diskrete Mathematik und Geomertie, Technische Universität Wien,
 nikolai.moshchevitin@tuwien.ac.at}
 \endgroup 
} 
 
\end{Large}
 
\vskip+1cm

In the present paper we give very simple general statements which deal with approximation of a real number by rationals and are related to isolation phenomenon.
In particular we study functions $ f_0(x)>f_1(x)>0$ such that existence of solutions
$\frac{p}{q}$
 of Diophantine inequality 
 $
\left| \alpha -\frac{p}{q}\right|< \frac{f_0(q)}{q^2}
$
leads to the existence of solutions of inequality
 $
\left| \alpha -\frac{p}{q}\right|< \frac{f_1(q)}{q^2}
$.
\vskip+0.1cm

In the present manuscript we correct a mistake from   our previous versions 
arXiv:2509.25628v1-3 and the published paper \cite{Pim}.

\vskip+0.1cm
In addition the new manuscript contains an Appendix with a result on exact simultaneous approximations to two numbers.

\section{Introduction}

Let $ \mathbb{Z}_+$ be the set of all positive integers and 
$$
\mathcal{Z}=
\{
m\in \mathbb{Z}_+:\,\,\, \text{exist}\,\, 
m_1,m_2 \in \mathbb{Z}_+\,\,\,
\text{such that }\,\,\, \max(m_1,m_2) \le m\,\,\,
\text{and}\,\,\,
m^2+m_1^2+m_2^2 = 3mm_1m_2
\} 
$$ 
be the set of  all Markoff numbers. We write elements of $\mathcal{Z}$ in increasing order as
$$
\mathcal{Z} :\,\,\,
m(1) = 1
<
m(2) =2
<
m (3)=5 
<...<
m(\nu)
<...
$$
and consider corresponding values
$$
\mu_1 = \frac{1}{\sqrt{5}}>
\mu_2  =\frac{1}{\sqrt{8}}>
\mu_3  =\frac{5}{\sqrt{221}}>...
>
\mu_\nu = 
\frac{1}{\sqrt{9- \frac{4}{m(\nu)^2}}}>
...
$$
One of the corollaries of the  famous theorem about the discrete part of Markoff and Lagrange spectra (see  \cite{cas},  Chapter II and \cite{Cus}) may be formulated as follows.

\vskip+0.3cm
{\bf Theorem A}. {\it
Let $ \gamma$ be a real number satisfying
$$
\mu_{\nu+1} < \gamma < \mu_\nu.
$$
If for a real number $\alpha$  inequality
 $$
\left| \alpha -\frac{p}{q}\right|< \frac{\gamma}{q^2}
$$
has infinitely many rational solutions $\frac{p}{q} \in \mathbb{Q}$,
then for this $\alpha$ a stronger inequality  
 $$
\left| \alpha -\frac{p}{q}\right|< \frac{\mu_{\nu+1} }{q^2}
$$
has infinitely many rational solutions $\frac{p}{q}$.
}
\vskip+0.3cm

One of the proofs of Theorem A contains as an ingredient a  special statement about "isolation" of the values of binary indefinite quadratic forms 
corresponding to Markoff numbers (see \cite{cas} Theorem I from Chapter II).

\vskip+0.3cm
On the other hand, quite recently Han\v{c}l \cite{H} proved a result which can be reformulated as follows.

\vskip+0.3cm
{\bf Theorem B}. {\it
For $x\ge 1$ consider function
$$
\psi_\nu (x) = \frac{2\mu_\nu}{\left(1+\sqrt{1+\frac{4\mu_\nu^2}{x^2}}\right) x^2}< \frac{\mu_\nu}{x^2}.
$$
If for a real number $\alpha$  inequality
 $$
\left| \alpha -\frac{p}{q}\right|< \frac{\mu_\nu}{q^2}
$$
has infinitely many rational solutions $\frac{p}{q} \in \mathbb{Q}$,
then for this $\alpha$ a stronger inequality  
 $$
\left| \alpha -\frac{p}{q}\right| \leqslant \psi_\nu (q)
$$
has infinitely many rational solutions $\frac{p}{q}$.
}
\vskip+0.3cm

The results mentioned above motivated us to prove a general statement concerning solutions of inequalities of the form
\begin{equation}\label{dsf}
\left| \alpha -\frac{p}{q}\right|< \frac{ f_0(q)}{q^2}
\end{equation}
with a positive function $f_0(x)<1$. Let us consider another function $f_1(x)$, satisfying
$$
 f_0(x) >f_1(x) >0, \,\,\, \text{for every}\,\,x.
$$
We are interested in the situation when there is no real $\alpha$ such that  simultaneously inequality 
(\ref{dsf})
has infinitely many rational solutions $\frac{p}{q}$ and inequality 
\begin{equation}\label{dsf1}
\left|\alpha -\frac{p}{q}\right|<\frac{f_1(q)}{q^2}
\end{equation}
 has no more than a finite number of solutions in rationals $\frac{p}{q}$;
in other words, existence of infinitely many solutions of
inequality (\ref{dsf}) leads to existence of infinitely many solutions of inequality (\ref{dsf1}). 
In our constructions function $f_0(x)$  can decrease to zero as quickly as desired.
These results may have relation to problems of mutual behaviour of two irrationality measure functions for real numbers 
(see \cite{KM, M, SR}).

\vskip+0.3cm
\section{Approximation to one real number}

First of all we formulate a general result for arbitrary functions $f(x)$.

\vskip+0.3cm
{\bf Theorem 1.}
  {\it Consider arbitrary function 
$f (x):\mathbb{Z}_+ \to (0,1]$. Then 
there exists $ \gamma \in(0,1]$  and a  function 
$g(x) :\mathbb{Z}_+ \to (0,1]$ such that 
$$
0< g(x) < \gamma \cdot f(x), \,\,\,\,\,
\text{for all}\,\,\, x\in \mathbb{Z}_+
$$ 
and
for any  real $\alpha$ the following statement holds:

\vskip+0.3cm
\noindent
if inequality 
 $$
\left| \alpha -\frac{p}{q}\right|< \frac{\gamma \cdot f(q)}{q^2}
$$
has infinitely many solutions in  rationals $\frac{p}{q} \in \mathbb{Q}$, then a stronger inequality
$$
\left| \alpha -\frac{p}{q}\right|< \frac{\gamma \cdot f(q)- g(q)}{q^2}
$$
also has infinitely many solutions   in  rationals $\frac{p}{q} \in \mathbb{Q}$.}
 
\vskip+0.3cm

Of course function $g(x)$ from Theorem 1 does not depend on $\alpha$. So we have a statement about solutions of inequalities (\ref{dsf}, \ref{dsf1}) with 

$$f_0(x) = \gamma f(x), \,\,\ f_1(x) = f_0(x) - g(x).$$

Theorem 1 is a particular case of a general statement which we formulate in Section \ref{fri} as Theorem 3.
To formulate this general result we need to discuss properties of Diophantine approximations for systems of linear 
forms in Section \ref{DAA} and consider some related functions in Section \ref{FF}. 
Theorem 1 will immediately follow from Theorem 3 and Lemma 2.

\vskip+0.3cm
In the case when $f(x)$ takes rational values only we have an effective version of Theorem 1.
Let 
 \begin{equation}\label{fq}
f(x) = \frac{A(x)}{B(x)}, \,\,\,\,
A(x)\in \mathbb{Z},\,\,\ \,B(x) \in \mathbb{Z}_+,\,\,\,\,
(A(x),B(x)) = 1.
\end{equation}
Define\footnote{ 
In our previous versions arXiv:2509.25628v1-3 the formula for $g_f(x)$ was incorrect, but this did not cause problems for the order of approximations; 
in particular, Example {\bf 1} remains just the same, in Example {\bf 2} the order of approximation remains the same and one needs only to change the constant, in Example { \bf 3}
there is a slight change in the exponent.}

 \begin{equation}\label{gq}
g_f(x)  = 
 \frac{(f(x))^2\cdot f\left( 3B(x) x^2\right)}{(A(x)x)^2}
 .
\end{equation}
It is clear that $ 0< g_f(x) <f(x)$ for every $x>1$. Moreover, if $f(x)$ is decreasing, then
 \begin{equation}\label{gqdec}
g_f(x)  \le
\frac{(f(x))^3}{ x^2}
 .
\end{equation}

\vskip+0.3cm
{\bf Theorem 2.} {\it
Let
 $f (x) : \mathbb{Z}_+ \to \mathbb Q_{+}\cap (0,1]$ be a non-strictly decreasing function.
Let $\alpha$ be real number. Assume that the inequality
\begin{equation}\label{4y}
\left| \alpha -\frac{p}{q}\right|< \frac{ f(q)}{q^2}
\end{equation}
has infinitely many solutions in  rationals $\frac{p}{q} \in \mathbb{Q}$. Then  for any positive $\varepsilon>0$ a stronger inequality
\begin{equation}\label{5y}
\left| \alpha -\frac{p}{q}\right|< \frac{ f(q)- (1-\varepsilon)g_f(q)}{q^2}
\end{equation}
where $g_f(x) $ is defined in (\ref{gq})
also has infinitely many solutions   in  rationals $\frac{p}{q} \in \mathbb{Q}$.}

\vskip+0.3cm

 We give a proof of Theorem 2 in Section \ref{t2p}.
 
\vskip+0.3cm
To illustrate Theorem 2 we consider three
examples.

\noindent
\textbf{ 1.} Let  $f(x)$  be just  a rational constant, that is for every $q$ we have  $f(q) = \frac{A}{B} \in \mathbb{Q}$.
Then
$g_f(x) = \frac{A}{B^3x^2}. $

\noindent
\textbf{ 2.} Let $f(x) = \frac{1}{x^{\sigma}},$ where $\sigma \in \mathbb{Z}_+.$
Then
$g_f(x) = \frac{1}{3^{\sigma} x^{\sigma^2 + 4 \sigma + 2}}.$

 \noindent
\textbf{ 3.} Let $f(x) = \frac{1}{2^x},$ then 
 $g_f(x) 
 =
 \frac{1}{x^2 2^{3x^2\cdot 2^x+2x}}
 $.

\vskip+0.3cm

{\bf Remark 1.}
 After the first version of  our paper was uploaded to arXiv system, S. Baker and B. Ward
 informed 
 us
 about their forthcoming manuscript \cite{bw} which deals with very related topics.
 In particular,
  a result similar to our Theorem 2 was also  obtained there. Moreover, \cite{bw} contains results on Hausdorff dimension of the related sets. The proofs seem to be rather different from ours.

\vskip+0.3cm

\section{Matrices and  approximations}\label{DAA}
Let $m,n\ge 1$ and 
\begin{equation}\label {vaa}
\Theta
=\left(
\begin{array}{ccc}
\theta_{1,1}&\cdots&\theta_{m,1}\cr
\cdots &\cdots &\cdots \cr
\theta_{1,n}&\cdots&\theta_{m,n}
\end{array}
\right)
\end{equation}
be a real matrix.
The main topic of the theory of linear Diophantine approximation is analysis of  behaviour of linear forms
$\sum_{i=1}^m\theta_{j}q_i - p_j$ when variables $q_i, 1\le i \le m $;      $p_j, 1\le j \le n$ take integer values.
We use notation
$$
q = (q_1,...,q_m) \in \mathbb{Z}^m,\,\,\,\,
p = (p_1,...,p_n) \in \mathbb{Z}^n,
\,\,\,\,
|q| = \max_{1\le i \le m} |q_i|,\,\,\,\,
||\Theta q|| = \min_{p\in \mathbb{Z}^n} |\Theta q - p|.
$$
The set of all matrices  (\ref{vaa}) with entries  $\theta_{i,j} \in [0,1]$ we denote by $\mathcal{R}$.
Fundamental Dirichlet's theorem claims that there exist infinitely many  $q\in \mathbb{Z}^m\setminus \{ (0,...,0)\}$ such that 
$$
||\Theta q|| \le  \frac{1}{|q|^{{m}/{n}}}.
$$
We deal with the set of positive integers $\mathbb{Z}_+$,
the set  of positive rationals $\mathbb{Q}_+$  and the set of positive reals $\mathbb{R}_+$.
Consider function
$f:\mathbb{Z}_+ \to (0,1]$.

\vskip+0.3cm

For integer  vector $q\in \mathbb{Z}^m\setminus \{ (0,...,0)\}$ define  sets
$$
M_q = M_q [f]=\left\{\Theta \in \mathcal{R}: \,\,  
||\Theta q||<  \frac{f(|q|)}{2 |q|^{{m}/{n}} }\right\},\,\,\,\,
{M}_q^c = M_q^c [f]= \mathcal{R}\setminus M_q
$$
and
$$
\mathcal{M} = \mathcal{M}[f]=  \bigcap_{Q=1}^\infty\,\,\, \bigcup_{q\in \mathbb{Z}^m: \, |q|\ge Q} M_q,\,\,\,\,\,\,
\mathcal{M}^c = \mathcal{M}^c[f]=   
\bigcup_{Q=1}^\infty
\mathcal{M}^c_Q,\,\,\,
\mathcal{M}^c_Q=
\mathcal{M}^c_Q[f]=
 \bigcap_{q\in \mathbb{Z}^m: \, |q|\ge Q} M_q^c.
$$

It is clear that $ \Theta \in \mathcal{M}$ if and only if 
\begin{equation}\label{1x}
\forall\, Q\,\,\, \exists \, q 
\in \mathbb{Z}^m
\,\,\,\text{such that}\,\,\,
|q|
\ge Q
\,\,\,
\text{and}
\,\,\,
 ||\Theta q || <  \frac{f(|q|)}{2|q|^{m/n}},
\end{equation}
meanwhile 
 $ \Theta \in \mathcal{M}^c$ if and only if 
\begin{equation}\label{2x}
\exists\, Q\,\,\, \forall \, q
\in \mathbb{Z}^m
\,\,\,\text{such that}\,\,\,
|q|
\ge Q\,\,\,\text{one has }\,\,\,
 ||\Theta q || \ge  \frac{f(|q|)}{2|q|^{m/n}}.
\end{equation}
In particular, as $f(x) >0$ for every $x \in \mathbb{Z}_+$, the set $\mathcal{M}^c$ does not contain any matrix $\Theta$ such that the {
 columns}
$$
\underbrace{
\left(
\begin{array}{c}
\theta_{1,1}
\cr
\vdots
\cr
\theta_{1,n}
\end{array}
\right),
\,\,\,
\dots
\,
,
\left(
\begin{array}{c}
\theta_{m,1}
\cr
\vdots
\cr
\theta_{m,n}
\end{array}
\right)}_m,
\,\,\,
\underbrace{
\left(
\begin{array}{c}
1
\cr
\vdots
\cr
0
\end{array}
\right),
\,\,\,
\dots
\,
,
\left(
\begin{array}{c}
0
\cr
\vdots
\cr
1
\end{array}
\right)}_n
$$
are linearly dependent over $\mathbb{Q}$.
 
\section{Family  
$
\frak{F}$}\label{FF}

For fixed $m,n$ we
consider a family of functions $f:\mathbb{Z}_+ \to (0,1]$ defined by
$$
\frak{F} = \{f:\mathbb{Z}_+ \to \mathbb{R}_+:\,\,\,\text{there exists}\,\,\, Q = Q[f]
\,\,\,
\text{such that for every}\,\,\, \Theta  \in  \mathcal{M}^c[f]
$$
$$\,\,\,\,\,\,\,\,\,\,\,\,\,\,\,\,\,\,\,\,\,\,\,\,
\text{and for every} \,\,\,
q\in \mathbb{Z}^m \,\,\,\text{with}\,\,\,  |q|\ge Q
\,\,\,
\text{one has }\,\,\,
|q| ^{m/n}\cdot ||\Theta q|| \neq f(|q|) \}.
$$
Here we would like to mention two examples of functions from the family 
$\frak{F}$.

\vskip+0.3cm
{\bf Lemma 1.}\footnote{ 
This Lemma is valid for $n=1$ only. More general formulation for arbitrary $n$ 
arXiv:2509.25628v1-3 and   \cite{Pim}
is incorrect, see Appendix.}

 Let $n = 1$ and function 
$ f(x): \mathbb{Z}_+\to (0,1] \cap \mathbb{Q} $. Then $f$ belongs to $\frak{F}$.
\vskip+1cm
Proof.

Let for some linear form \(\Theta\) and \(q_{0} \in \mathbb{Z}^m\) we have
$
\| \Theta q_{0}\| = \frac{f(|q_0|)}{|q_0|^m}.
$
Since the function \(f\) takes only rational values, we get
$
\frac{f(|q_0|)}{|q_0|^m} \in \mathbb{Q},
$
hence \(\Theta q_{0} \in \mathbb{Q}\). Let
$$
\Theta q_0 = \frac{A}{B}, \qquad (A, B) = 1.
$$
Then for every positive integer \(k\) we have 
$$
\| \Theta \cdot q_{0}Bk\| = \|Ak \| = 0
< \frac{f(|q_{0}Bk|)}{2|q_{0}Bk|^{m}}.
$$

Thus, there exist infinitely many vectors \(q\) satisfying the condition
$\| \Theta q\| < \frac{f(|q|)}{2|q|^{m}},$
and hence \(\Theta \notin \mathcal{M}^c[f]\).

Thus, \(f \in \mathfrak{F}\), as required.
$\Box$
 
 \vskip+0.3cm
Our next lemma deals with a special case $m=n=1$.
\vskip+0.3cm
{\bf Lemma 2.} Let $m=n=1$. Then for any function  $f:\mathbb{Z}_+ \to (0,1] $ there exists  
$\gamma \in  \left[\frac{1}{2},1\right]$
such that $\gamma \cdot f \in \frak{F}.$

\vskip+0.3cm

Before giving further proofs we need to introduce auxiliary sets
$$
V_{q,p} = \left\{ \alpha \in [0,1]: \,\,  
 |\alpha q - p|< \frac{f(q)}{q} \right\},\,\,\,\,\,
W_{q,p} = \left\{ \alpha \in [0,1]: \,\,  
 |\alpha q - p|<\frac{f(q)}{4q} \right\},
$$
 and the sets
$$
\mathcal{K} = \mathcal{K}[f]=  \bigcap_{Q=1}^\infty
\bigcup_{q\ge Q}  \bigcup_{p\in \mathbb{Z}} V_{q,p}
,\,\,\,\,\,\,
\mathcal{N}_Q  =
\mathcal{N}_Q [f] =
[0;1]
\setminus
\left(\bigcup_{ q\ge Q}\,\,\, 
\bigcup_{p\in \mathbb{Z}} 
 W_{q,p}  \right)
  = \bigcap_{q\ge Q} 
 \left\{\alpha \in [0;1]: \,\,  
||\alpha q||\ge  \frac{f(q)}{ 4q}\right\}.
$$
We should note  that for every $Q$  the set $
\mathcal{N}_Q$
is non-where dense in $[0;1]$. 
 
Also   for every $\gamma \in \left[\frac{1}{2},1\right]$ we have
\begin{equation}\label{halbe1}
\mathcal{M}[\gamma \cdot f]  
\subset
 \mathcal{K}[f] 
\end{equation}
and
\begin{equation}\label{halbe}
\mathcal{M}^c_Q[\gamma \cdot f]  \subset \mathcal{N}_Q[f].
\end{equation}

\vskip+0.3cm
 \vskip+0.3cm
Proof of Lemma 2.   Note that  for $ \gamma \in \left[\frac{1}{2},1\right]$ condition 
$\gamma \cdot f \in \frak{F}$  holds if and only if
 \begin{equation}\label{ifandonly}
\gamma \in 
\bigcup_{Q=1}^\infty
\bigcap_{q: |q|\ge Q} {B}_{q},
\,\,\,
\text{where}
\,\,\,
 {B}_{q} =
\left\{\beta \in  \left[\frac{1}{2},1\right]:   
\beta \neq \frac{q}{ 
f(q)
}
\cdot ||\alpha q||
\,\,\,
\forall \, \alpha \in \mathcal{M}^c[\beta\cdot f] \right\}.
\end{equation}
For integers $Q,q$ define
 $$
\mathcal{H}_{Q;q} =
 \left\{\beta \in  \left[\frac{1}{2},1\right]:   
 \exists\, \alpha \in \mathcal{M}^c_Q[\beta\cdot f] \,\,
 \text{such that}\,\,
\beta = \frac{q}{ 
f(q)
}
\cdot ||\alpha q||
 \right\}
.
$$
Then
$$
B_q =    \left[\frac{1}{2},1\right] \setminus 
\left(\bigcup_{Q=1}^\infty
\mathcal{H}_{Q;q}
\right).
 $$
So if
 $ \gamma \in \left[\frac{1}{2},1\right]$ satisfies
$\gamma \cdot f \not\in \frak{F}$   
then
$$
\gamma \in \bigcap_{Q=1}^\infty\bigcup_{q\ge Q}  
\left(\left[\frac{1}{2},1\right] \setminus B_q\right)
\subset
 \bigcup_{Q=1}^\infty \bigcup_{q=1}^\infty  \mathcal{H}_{Q;q}
  \subset
  \bigcup_{Q=1}^\infty\bigcup_{q=1}^\infty \left\{\beta \in  \left[\frac{1}{2},1\right]:   
\exists\, \alpha \in \mathcal{N}_Q[ f] 
\,\,\text{such that}\,\,
\beta=
 \frac{q}{ 
f(q)} ||\alpha q||
\right\}
$$
(the last inclusion here is due to 
 (\ref{halbe})).
For integers $Q,q,p$ consider the set
$$
\mathcal{N}_{Q;q,p} = \left\{
\beta \in  \left[\frac{1}{2},1\right]:   
\exists\, \alpha \in \mathcal{N}_Q[ f] 
\,\,\text{such that}\,\,
\beta = \frac{q}{f(q)} \cdot |q\alpha - p|\right\}.
$$
This set is also non-where dense.
We see that  if
 $ \gamma \in \left[\frac{1}{2},1\right]$ satisfies
$\gamma \cdot f \not\in \frak{F}$   
then

\begin{equation}\label{3x}
\gamma \in \bigcup_Q\bigcup_{q}\bigcup_{p} \mathcal{N}_{Q;q,p}.
\end{equation}
The right hand side of (\ref{3x}) is a countable union of non-where dense sets. By
Baire's theorem there exist $\gamma \in \left[\frac{1}{2},1\right]$  which does not belong to this set.$\Box$

\vskip+0.3cm

{\bf Remark 2.}
From the proof of Lemma 2 we see that the set of all $\gamma$ satisfying the conclusion of Lemma 2 (and thus satisfying the conclusion of Theorem 1) is a set of second Baire category. 

\vskip+0.3cm
 
{\bf Remark 3.}
Let  
$\lim_{x\to \infty} f(x) = 0$ and 
the series  
\begin{equation}\label{halbe2}
\sum_{x=1}^\infty \frac{f(x)}{x}
\end{equation}
 converges. Then for almost all 
$\gamma$ in the sense of Lebesgue measure one has $\gamma \cdot f \not\in \frak{F}$.

\vskip+0.3cm
Proof of Remark 3.
It is enough to consider the case $\gamma \in \left[\frac{1}{2},1\right]$.
 As   for $\gamma$  under condition $\gamma \cdot f \in \frak{F}$  we have (\ref{ifandonly}),
  to prove {Remark 3} it is enough to show  for $q$ large enough  each ${B}_q$ is a set of zero measure.

For non-zero integer $q \in\mathbb{Z}$  consider the map
$$
T_q: [0,1]\to\mathbb{R},\,\,\,
T_q(\alpha) = \frac{q}{f(q)} ||\alpha q||.
$$
As $f(x)$ tends to zero, for $q$ large enough  the image  $T_q ([0,1])$ contains the interval 
$\left[\frac{1}{2},1\right]$. So for $\beta \in \left[\frac{1}{2},1\right]$ we have
$$
\beta \neq \frac{q}{ 
f(q)
}
\cdot ||\alpha q||
\,\,\,
\forall \, \alpha \in \mathcal{M}^c[\beta\cdot f]
\,\,\,
\Longleftrightarrow
\beta \not\in  T_q  \left(\mathcal{M}^c[\beta\cdot f]\right)
\,\,\,
\Longrightarrow
\,\,\,
\beta \in T_q  \left(\mathcal{M}[\beta\cdot f]\right).
$$

By (\ref{halbe1}) we see that  for any $q$  the inclusion
$
B_q \subset
T_q  \left(\mathcal{M}[\beta\cdot f]\right)\subset T_q  \left(\mathcal{K}[ f]\right)
 $
 hold. But we should note that   
if the series
(\ref{halbe2})
converges, the  set $ \mathcal{K}[f]$ has zero Lebesgue measure,
by standard Borel-Cantelli argument (see \cite{Spri}).
 So  for every $q$ large enough $
B_q$ is also a zero set.$\Box$
\section{General result and its proof}\label{fri}

{\bf Theorem 3.} {\it Let $f \in \frak{F}$. Then there exists a positive valued function
$g:\mathbb{Z}_+\to (0,1]$ such that for any  real matrix $\Theta $ the following statement is valid:

if inequality 
\begin{equation}\label{4x}
|| \Theta q||< \frac{f(|q|)}{|q|^{m/n}}
\end{equation}
has infinitely many solutions in  ${q}\in \mathbb{Z}^m$ then a stronger inequality
\begin{equation}\label{5x}
||\Theta q ||< \frac{f(|q|)- g(|q|)}{|q|^{m/n}}
\end{equation}
also has infinitely many solutions  in ${q}\in \mathbb{Z}^m$.}

\vskip+0.3cm
Lemma 1 from the previous section leads to the following

\vskip+0.3cm
{\bf Corollary.}\footnote{
This Corollary  is valid for $n=1$ only. More general formulation for arbitrary $n$  from 
arXiv:2509.25628v1-3 and   \cite{Pim}
is incorrect, see Appendix.}
 {\it  Let $n =1$.
Let \(f(x):\mathbb{Z}_+ \to (0,1] \cap \mathbb{Q}\). Then there exists a function \(g:\mathbb{Z}_+ \to (0,1]\) such that for any real linear form \(\Theta\), the existence of infinitely many solutions of the inequality
\[
\| \Theta q\| < \frac{f(|q|)}{|q|^{m}}
\]
implies the existence of infinitely many solutions of the stronger inequality
\[
\| \Theta q\| < \frac{f(|q|)-g(|q|)}{|q|^{m}}.
\]. 
 }

\vskip+0.3cm
Lemma 2 shows that Theorem 1 formulated in the beginning of our paper immediately follows from Theorem 3.

\vskip+0.3cm
Proof of Theorem 3. 
Define $g_1(x)  = \frac{f(x)}{2}$.
Each set $ \mathcal{M}_x^c \subset \mathcal{R}$ is a compact set. As $f\in \frak{F}$, for $ x\ge Q[f]$ we  define 
\begin{equation}\label{7x}
  g_2(x) = \frac{1}{2} \cdot \min_{q\in \mathbb{Z}^m: |q|= x}\,\, \min_{\Theta \in \mathcal{M}_{x}^{c}}
    \left| f(x) - x^{m/n}\cdot ||\Theta q ||\right| > 0
\end{equation}
The minimum is attained because each set \(\mathcal{M}_{x}^{c}\) is compact, and the positivity of the quantity follows from the fact that the function \(f\) belongs to the set \(\mathfrak{F}\).

Again we consider two cases. In the first case we assume that 
$\Theta\in \mathcal{M}$. Then we have the result of Theorem 2 with $ g(q) = g_1(q) $ because of 
(\ref{1x}).

Consider the case $\Theta\in \mathcal{M}^c$.

Then there exists an infinite sequence of integer vectors $q$ satisfying (\ref{4x}) 
and such that
  $ \Theta \in \mathcal{M}^c_{|q|}$ (see (\ref{2x})). By
(\ref{7x}) for these values of $q$ with $ |q|\ge Q[f]$ we have
$$ |q|^{m/n}\cdot ||\Theta q|| \le  f(|q|) - 2g_2(|q|) <f(|q|) - g_2(|q|).
$$

So we have the conclusion of Theorem 3 with   $ g(x) = \min (g_1(x), g_2(x))$.$\Box$

\section{Proof of Theorem 2}\label{t2p}

We may assume that $\alpha$ is irrational. Consider its representation as an ordinary  infinite continued fraction 
$$
\alpha = [a_0;a_1,a_2,...,a_n,...],\,\,\,\,\, a_0 \in \mathbb{Z}, \,\,\, a_j \in \mathbb{Z}_+, j \ge 1
$$
and define
$$
\frac{p_n}{q_n} =  [a_0;a_1,a_2,...,a_n], \,\,\,\,\,
\alpha_{n+1} = [a_{n+1};a_{n+2},...,a_{n+k},...].
$$
It is well known  (see for example \cite{RS}, Chapter 1, \S 6)  that 
\begin{equation}\label{perron0}
\alpha_n^* := \frac{q_{n-1}}{q_n} = [0;a_n,a_{n-1},...,a_{1}]
\end{equation}
and
\begin{equation}\label{perron}
|q_n \alpha - p_n| = \frac{1}{q_n (\alpha_{n+1}+\alpha_n^*)}.
\end{equation}
Define
$$
\frac{p_{n,k}}{q_{n,k}} =  [a_{n+1};a_{n+2},a_{n+3},...,a_{n+k}],\,\,\,\,\,
\alpha_{n,k}^* = \frac{q_{n,k-1}}{q_{n,k}}= [0;a_{n+k},a_{n+k-1},....,a_{n+2}]
.
$$
Then application of (\ref{perron}) gives
\begin{equation}\label{perron1}
|q_{n,k} \alpha_{n+1} - p_{n,k}| = \frac{1}{q_{n,k} (\alpha_{n+k+1}+\alpha_{n,k}^*)}.
\end{equation}

 \vskip+0.3cm
 Let us denote $ f_n = f(q_n)$.
We will distinguish two cases.
In the {\bf first case}
we assume that there exist infinitely many values of $n$ such that 
\begin{equation}\label{c1}
\frac{1}{\alpha_{n+1}+\alpha_n^*} < f_n \left(1-\left(\frac{f_n}{q_n}\right)^2\right).
\end{equation}
We see that in the {\bf first case} condition (\ref{c1}) together with (\ref{gqdec})  and (\ref{perron}) immediately gives infinitely many solutions of inequality (\ref{5y}).

 \vskip+0.3cm
 Now we consider the   {\bf second case}, when inequality (\ref{c1}) is valid for not more than for a finite collection of values of $n$, Then there exists $n_0$ such that for all 
 $n\ge n_0$ one has
 \begin{equation}\label{c2}
\frac{1}{\alpha_{n+1}+\alpha_n^*} \ge f_n \left(1-\left(\frac{f_n}{q_n}\right)^2\right).
\end{equation}
In particular,
 \begin{equation}\label{c22}
 \frac{q_{n+1}}{q_n} = a_{n+1}+\alpha_n^* < \alpha_{n+1}+\alpha_n^*\le \frac{1}{f_n \left(1-\left(\frac{f_n}{q_n}\right)^2\right)}
 .
\end{equation}

By condition  (\ref{4y})  and formula (\ref{perron}) we may assume that  for infinitely many values of $n\ge n_0$ we have 
 \begin{equation}\label{c3}
\frac{1}{\alpha_{n+1}+\alpha_n^*} <f_n.
\end{equation}
In the rest of the proof we deal just with these values of $n$.
We consider the value 
$$
0<
\delta_n= 
\alpha_{n+1}+\alpha_n^* -\frac{1}{f_n}
=
\alpha_{n+1}  +\frac{q_{n-1}}{q_n} - \frac{1}{f_n} =
\alpha_{n+1} - \frac{C_{n}}{A_nq_n},
$$
{where}
$$
 C_n = 
 q_{n}B_{n}  - q_{n-1}A_{n} \in \mathbb{Z},\,\,\, A_n = A(q_n),\,\,\, B_n = B(q_n)
 $$
 (we used here (\ref{fq}) and (\ref{perron0})).
 
 Let us choose $k$ such that 
  \begin{equation}\label{c31}
 q_{n,k} 
 \le A_nq_n < q_{n,k+1}.
 \end{equation}
 Then by Lagrange's theorem on best approximations   we get
 \begin{equation}\label{c5}
A_nq_n\delta_n \ge |q_{n,k}\alpha_{n+1} -p_{n,k}|.
\end{equation}
Now  from inequality (\ref{c5}), equality  (\ref{perron1}), inequalities (\ref{c2}) and
$$
|\alpha_{n+k}^* - \alpha_{n,k}^*|\le \frac{1}{q_{n,k}^2},
$$
we get
$$
A_nq_n\delta_n \ge
\frac{1}{q_{n,k}(\alpha_{n+k+1} + \alpha_{n,k}^*)}\ge
\frac{1}{q_{n,k}(\alpha_{n+k+1} + \alpha_{n+k}^*+1/q_{n,k}^2)}\ge
\frac{1}{
 \frac{q_{n,k}  }{f_{n+k}\left( 1- \left(\frac{f_{n+k}}{q_{n+k}}\right)^2\right)}+\frac{1}{q_{n,k}}}
 =
 $$
 $$
  \frac{f_{n+k}}{q_{n,k}} \left(1- \Big(\frac{f_{n+k}}{ q_{n+k}}\Big)^2\right) \frac{1}{ 1+ \frac{f_{n+k} } {q_{n,k}^2 } \left( 1- \left(\frac{f_{n+k}}{q_{n+k}}\right)^2\right)}
> \frac{f_{n+k}}{q_{n,k}} \left(1-\frac{2f_{n+k}}{ q_{n,k}^2}\right)
$$
for sufficiently large \(n\).
Then
\begin{equation}\label{continuant}
q_{n+k} \le (q_{n+1}+q_n) q_{n,k}
\le q_n q_{n,k} \left( \frac{q_{n+1}}{q_n}+1\right) \le
q_n q_{n,k} \cdot \left(\frac{1}{f_n \left(1-\left(\frac{f_n}{q_n}\right)^2\right)}+1\right),
\end{equation}
by (\ref{c22}). Also, for \(q_{n} \ge 2\) we have

\begin{equation}\label{kazah}
    A_{n}q_{n}^2 \Big(\frac{1}{f_{n}(1 - (\frac{f_n}{q_n})^{2})}+1\Big) \le 3B_{n}q_{n}^2.
\end{equation}

Now, taking into account (\ref{continuant}), (\ref{kazah}) and the monotonicity of \(f(x)\), we finally obtain for sufficiently large \(n\)
$$
\delta_n \ge \frac{f\left(A_n q_n^2\cdot \left(\frac{1}{f_n \left(1-\left(\frac{f_n}{q_n}\right)^2\right)}+1\right) \right)}{(A_nq_n)^2} \cdot \left(1-\frac{2f_{n+k}}{ q_{n,k}^2}\right) > \frac{f(3B_n q_{n}^2)}{(A_{n}q_{n})^2} \cdot R_n,
$$
where

$$R_n = \frac{A_{n}q_{n}}{q_{n,k}}\Big(1 - \frac{2f_{n+k}}{q_{n,k}^2} \Big).$$

Let \(\varepsilon \in (0;1)\) and let a positive integer \(M\) be such that \(\frac{2}{M^2} \le \varepsilon.\) Consider the cases.

\begin{itemize}
    \item If \(q_{n,k} \ge M,\) then \(R_n \ge 1 - \frac{2}{M^2} \ge 1 - \varepsilon.\)

    \item If \(2 \le q_{n,k} \le M,\) then \(R_n \ge \frac{A_{n}q_{n}}{2M} > 1\) for \(q_n > 2M.\)

    \item If \(q_{n,k} = 1,\) then from (\ref{c22}) we have

    $$a_{n+k+1} < \frac{q_{n+k+1}}{q_{n+k}} \le \frac{1}{f_{n+k}(1 - (\frac{f_{n+k}}{q_{n+k}})^2)} \le \frac{2}{f_{n+k}},$$

    and using (\ref{c31}) we get \(R_n = A_{n}q_{n}(1- 2f_{n+k}) \ge q_{n} - 4 \ge 1\), when \(q_{n} \ge 5 \).
\end{itemize}

Thus, for sufficiently large \(n\)

$$\delta_n > (1-\varepsilon) \frac{f\left(3B_{n}q_{n}^2 \right)}{(A_nq_n)^2}.$$

Since

$$\alpha_{n+1}+\alpha_n^* =\frac{1}{f_n}+\delta_n$$

using (\ref{perron}), we conclude that for any \(\varepsilon>0\), for the large \(n\) considered by us, we have
$$
\left|\alpha - \frac{p_n}{q_n}\right| = \frac{f_{n}}{q_{n}^2} \Big( 1 - \frac{f_{n} \delta_n}{1 + f_{n} \delta_n} \Big)< \frac{f_n}{q_n^2} \left( 1- (1-\varepsilon) \frac{ f_n \cdot f\left(3B_{n}q_{n}^2 \right)}{(A_nq_n)^2}\right).
$$
Theorem is proven.$\Box$

  \section{Appendix}
  
  In this section we will show that Corollary to Theorem 3 cannot be generalised to the case $n\ge 2$.
  Moreover, to show this  we give here the final result about existence of two 
  real numbers admitting  exact simultaneous approximations under a condition on of the type $ f(x) = O(x^{-1})$.  
  For such type of functions the result is essentially better than that from paper  \cite{A}.
  
    \vskip+0.3cm
  
For real $\alpha,\beta$ and a positive integer $q$ we write
$\Theta = \left(\begin{array}{c}\alpha\cr\beta\end{array}\right) $ and $
\|q\Theta \|=\max\big(\|q\alpha\|,\|q\beta\|\big),
$
where for a real number  $\|\cdot\|$ is  just the ordinary the distance to the nearest integer.
 
   \vskip+0.3cm
   {\bf Theorem A1.}
 {\it 
Let $f:\mathbb{Z}_+\to(0,1]$ be a non-increasing function such that
\begin{equation}\label{f}
f(q)\le\frac1{3q}\qquad\text{for all }q\ge q_0 ,
\end{equation}
and let $g:\mathbb{Z}_+\to\mathbb{R}$ satisfy $0<g(q)<f(q)$ for all $q$. Then there exist irrational
numbers $\alpha,\beta$ and an infinite sequence $\frak{Q}$ of positive integers $q_1<q_2<q_3< \cdots$,
$q_1\ge q_0$, such that
\begin{equation}
\|q_k\Theta \|<f(q_k)\qquad\text{for all }k,\label{A}
\end{equation} whereas 
\begin{equation}\label{AW}
 \text{
 inequality}\,\,\, \|q\Theta \|<f(q)-g(q)\,\,\, \text{ has no solutions} \,\,\,
q\ge q_1.
\end{equation}
}
  \vskip+0.3cm
  
  In the proof of Theorem A1 to ensure (\ref{AW}) we will show that 
\begin{equation}
\|q\Theta \|\ge f(q)\qquad\text{for all }q\ge q_1,\ q\notin\frak{Q},\label{B}
\end{equation}
and
\begin{equation}
\|q_k\Theta \|\ge f(q_k)-g(q_k)\qquad\text{for all }k.\label{C}
\end{equation}
 
\vskip+0.3cm
Proof of Theorem A1.
 
Let $\beta\in(0,1)$ be irrational, $\beta=[0;b_1,b_2,b_3,\dots]$, with convergents
$P_\nu/Q_\nu$, $Q_0=1$, $Q_1=b_1$, $Q_{\nu+1}=b_{\nu+1}Q_\nu+Q_{\nu-1}$, and
$\beta_{\nu+1}=[b_{\nu+1};b_{\nu+2},\dots]$.From continued fractions theory we have\begin{equation}\label{cf}
\|Q_\nu\beta\|=|Q_\nu\beta-P_\nu|=\frac1{Q_\nu\beta_{\nu+1}+Q_{\nu-1}},\qquad\text{hence}\qquad
\frac1{Q_\nu(b_{\nu+1}+2)}<\|Q_\nu\beta\|<\frac1{Q_\nu b_{\nu+1}} ,\,\,\,\, \nu \ge 1.
\end{equation}
We formulate two easy auxiliary statements. The first one immediately follows from Legendre Theorem. The second one is obvious.
 
   \vskip+0.3cm
   {\bf Lemma A2.} {\it
Let $q\ge1$ and let $a$ be the integer nearest to $q\beta$. If
$(q,a)\ne\lambda\,(Q_\nu,P_\nu)$ for all $\nu\ge1$ and all $\lambda\in\mathbb{Z}_+$, then
$
\|q\beta\|\ \ge\ \frac1{2q}.
$
If $(q,a)=\lambda\,(Q_\nu,P_\nu)$, then $\|q\beta\|=\lambda\|Q_\nu\beta\|$.
}
 
   \vskip+0.3cm
   {\bf Lemma A3.} {\it
Let $\alpha$ be real, $q\ge2$, $a\in\mathbb{Z}$, $\varepsilon\in(0,1)$, $\xi\in(0,\tfrac12)$, and
\begin{equation}\label{window}
\alpha\in\Big[\frac aq+\frac{(1-\varepsilon)\xi}{q},\ \frac aq+\frac{(1-\varepsilon/2)\xi}{q}\Big].
\end{equation}
Then $(1-\varepsilon)\xi\le\|q\alpha\|\le(1-\varepsilon/2)\xi$.
}

   \vskip+0.3cm
 
 Now we describe the construction of numbers $\beta$ and $\alpha$.
 
    \vskip+0.3cm
Put
\begin{equation}\label{eps}
\varepsilon_k=\frac{g(q_k)}{f(q_k)}\in(0,1),\qquad
b_k=\Big\lceil\frac1{q_kf(q_k)}\Big\rceil ,
\end{equation}
where the integers $q_k$ are defined below. The number $\beta$ is
\begin{equation}\label{beta}
\beta=[\,0;\ \underbrace{1,\dots,1}_{n_1},\frak{b}_1,\ \underbrace{1,\dots,1}_{n_2},\,\frak{b}_2,\ \dots\,],
\qquad n_k\ge1,
\end{equation}
so that $\frak{b}_k=b_{\nu_k+1}$ with
\[
\nu_k=n_1+\dots+n_k+k-1,\qquad q_k:=Q_{\nu_k} .
\]
The number $\alpha$ will be  the common point of nested intervals
\begin{equation}\label{Ik}
I_k=\Big[\frac{a_k}{q_k}+\frac{(1-\varepsilon_k)f(q_k)}{q_k},\
\frac{a_k}{q_k}+\frac{(1-\varepsilon_k/2)f(q_k)}{q_k}\Big],\qquad
|I_k|=\frac{\varepsilon_kf(q_k)}{2q_k} 
\end{equation}
constructed inductively.

Now we describe the inductive procedure to choose
integers $n_k$, $a_k$. Let $n_1$ be the least integer with
$Q_{n_1}\ge\max(q_0,3)$ (here $Q_{n_1}$ is a Fibonacci number), $q_1=Q_{n_1}$, $a_1=1$;
this defines $\frak{b}_1$, $\varepsilon_1$, $I_1\subset(0,1)$. Suppose $n_1,\dots,n_k$ and
$a_1,\dots,a_k$ are chosen, so that $b_1,\dots,b_{\nu_k+1}$ and $I_k$ are determined. Let
$n_{k+1}$ be the least integer $n\ge1$ such that the denominator $Q_{\nu_k+1+n}$ of the
continued fraction $[0;b_1,\dots,b_{\nu_k+1},\underbrace{1,\dots,1}_{n}]$ satisfies
\begin{equation}\label{big}
Q_{\nu_k+1+n}\ \ge\ \frac3{|I_k|} .
\end{equation}
Put $q_{k+1}=Q_{\nu_{k+1}}$. By \eqref{big} the points $a/q_{k+1}$, $a\in\mathbb{Z}$, form a
progression with step $\le|I_k|/3$, so the middle third of $I_k$ contains such a point; let
$a_{k+1}$ be the least admissible $a$. Then \eqref{Ik} defines $I_{k+1}$, and since
$f\le1$ and $1/q_{k+1}\le|I_k|/3$,
\[
I_{k+1}\subset\Big[\frac{a_{k+1}}{q_{k+1}},\frac{a_{k+1}+1}{q_{k+1}}\Big]\subset I_k .
\]
Finally $\alpha$ is the unique point of $\bigcap_kI_k$.
 
The number $\beta$ is irrational, having an infinite continued fraction. The number
$\alpha$ is irrational: the fractions $a_k/q_k$ are pairwise distinct   and $|\alpha-a_k/q_k|<f(q_k)/q_k\le1/(3q_k^2)$, which is
impossible for a rational $\alpha$ when $q_k\to\infty$. Moreover, by the construction suggested one can easily  ensure
$1,\alpha, \beta$ to be linearly independent over $\mathbb{Q}$.

 \vskip+0.3cm
Finally, we verify the conclusion of Theorem A1.
 
\emph{Inequalities \eqref{A} and \eqref{C}.} By \eqref{cf} and \eqref{eps},
$\|q_k\beta\|<1/(q_k\frak{b}_k)\le f(q_k)$. By Lemma A3 applied with
$\xi=f(q_k)$, $\varepsilon=\varepsilon_k$ (note $\alpha\in I_k$),
\[
f(q_k)-g(q_k)=(1-\varepsilon_k)f(q_k)\le\|q_k\alpha\|\le(1-\varepsilon_k/2)f(q_k)<f(q_k).
\]
Hence $f(q_k)-g(q_k)\le\|q_k\Theta\|<f(q_k)$.
 
\emph{Inequality \eqref{B}.} Let $q\ge q_1$, $q\notin\frak{Q}$, and let $a$ be the integer
nearest to $q\beta$. We show $\|q\beta\|\ge f(q)$; the number $\alpha$ is not needed.
Recall $f(q)\le1/(3q)$ by \eqref{f}.
 
1°. $(q,a)$ is not proportional to a convergent pair of $\beta$. Then by Lemma A2.
\[
\|q\beta\|\ge\frac1{2q}>\frac1{3q}\ge f(q).
\]
 
2°. $(q,a)=\lambda(Q_\nu,P_\nu)$ with $\nu\ge1$, $\lambda\ge1$, so that
$\|q\beta\|=\lambda\|Q_\nu\beta\|$ by Lemma A2.
 
(a) $\nu\notin\{\nu_k\}$. Then $b_{\nu+1}=1$, and \eqref{cf} gives
\[
\|q\beta\|=\lambda\|Q_\nu\beta\|>\frac{\lambda}{3Q_\nu}=\frac{\lambda^2}{3q}\ge\frac1{3q}\ge f(q).
\]
 
(b) $\nu=\nu_k$ for some $k$, $\lambda\ge2$. Then $Q_\nu=q_k$, $b_{\nu+1}=\frak{b}_k$, and by
\eqref{cf}, \eqref{eps} and $3q_kf(q_k)\le1$,
\[
\|q_k\beta\|>\frac1{q_k(\frak{b}_k+2)}\ge\frac1{q_k\big(\frac1{q_kf(q_k)}+3\big)}
=\frac{f(q_k)}{1+3q_kf(q_k)}\ge\frac{f(q_k)}2 .
\]
Hence, as $q=\lambda q_k\ge q_k$ and $f$ is non-increasing,
\[
\|q\beta\|=\lambda\|q_k\beta\|>\frac{\lambda f(q_k)}{2}\ge f(q_k)\ge f(q).
\]
 
(c) $\nu=\nu_k$, $\lambda=1$. Then $q=q_k$, which is excluded.
 
Thus $\|q\Theta\|\ge\|q\beta\|\ge f(q)$ in all cases, and \eqref{B} is proved.
Theorem A1 is proved. \qed

  \vskip+0.3cm

  \noindent
  {\bf Acknowledgements}
  \,\, Sergei Pitcyn's research is supported by
  Russian Scientific Foundation, grant
  no.~25-11-00112, https://rscf.ru/project/25-11-00112/.
 Nikolay Moshchevitin’s research is supported by Austrian Science Fund (FWF), Forschungsprojekt PAT1961524.
 
   \vskip+0.3cm

  \noindent
  {\bf AI declaration.} During preparation of this  corrected version of the manuscript 
  the authors used Claude Anthropic, version Fable 5.1.

\end{document}